\documentclass[a4paper,fleqn]{article}

\usepackage{fullpage}

\usepackage{times,cite,w-thm}
\usepackage[T1]{fontenc}
\usepackage[utf8]{inputenc}
\usepackage{graphicx}
\usepackage{amsthm}
\usepackage{amsmath,amsfonts}
\def\Sm{C_\Theta^{\infty}}
\def\D{D(\mathcal{I})}

\title{A new matrix equation expression for the solution of non-autonomous linear systems of ODEs}
\providecommand{\keywords}[1]{\textbf{\hspace{0.9cm}Keywords: } #1}
\date{}
\author{ Stefano Pozza \footnotemark[2]	\and Niel Van Buggenhout\footnotemark[2]}

\begin{document}
	\maketitle
		\renewcommand{\thefootnote}{\fnsymbol{footnote}}
	
	\footnotetext[2]{Charles University, Sokolovská 83 186, 75 Praha 8, Czech Republic. (pozza@karlin.mff.cuni.cz, buggenhout@karlin.mff.cuni.cz )}
	\footnotetext{This work was supported by Charles University Research programs No. PRIMUS/21/SCI/009 and UNCE/SCI/023, and by the Magica project ANR-20-CE29-0007 funded by the French National Research Agency.}

\begin{abstract}
The solution of systems of non-autonomous linear ordinary differential equations is crucial in a variety of applications, such us nuclear magnetic resonance spectroscopy. A new method with spectral accuracy has been recently introduced in the scalar case. The method is based on a product that generalizes the convolution. In this work, we show that it is possible to extend the method to solve systems of non-autonomous linear ordinary differential equations (ODEs). In this new approach, the ODE solution can be expressed through a linear system that can be equivalently rewritten as a matrix equation. Numerical examples illustrate the method's efficacy and the low-rank property of the matrix equation solution.

\end{abstract}
\maketitle                   

\section{Introduction}
Systems of non-autonomous linear ordinary differential equations (ODEs) appear in a variety of applications, 
and its numerical computation is often challenging, particularly for large-to-huge size systems.
	For instance, in nuclear magnetic resonance spectroscopy (NMR) \cite{HaSp98}, the system solution describes the dynamics of the nuclear spins of a sample in a time-varying magnetic field. 
	The size of such systems is $2^k \times 2^k$ for a sample with $k$ spins and is usually sparse.
In \cite{PozVan_PANM21}, we proposed a new method with spectral accuracy for solving \emph{scalar} non-autonomous ordinary differential equations.
 In the present work, we extend this method to the case of systems of non-autonomous ODEs.

Consider a matrix $\tilde{A}(t) \in \mathbb{C}^{N \times N}$ composed of elements from $C^\infty(\mathcal{I})$, i.e., the set of functions infinitely differentiable (smooth) over $\mathcal{I}$, with $\mathcal{I}$ a closed and bounded interval in $\mathbb{R}$. The system
\begin{equation}\label{eq:ode:intro} 
\frac{d}{dt}U_s(t) = \tilde{A}(t) U_s(t), \quad U_s(s)=I_N, \quad \text{ for } t \geq s, \quad t,s\in \mathcal{I},
\end{equation}
has a unique solution $U_s(t) \in \mathbb{C}^{N \times N}$; $I_N$ stands for the $N \times N$ identity matrix.
Note that the condition $U_s(s) = I_N$ is not restrictive, since, given a matrix $B \in \mathbb{C}^{N\times N}$, the matrix-valued function $V_s(t) := U_s(t)B$ solves the ODE 
\begin{equation*}
  \frac{d}{dt} V_s(t) = \tilde{A}(t) V_s(t), \quad V_s(s) = B \quad \text{ for } t \geq s, \quad t,s\in \mathcal{I}. 
\end{equation*}

At the heart of the new method for solving \eqref{eq:ode:intro} is a non-commutative convolution-like product, denoted by $\star$, defined between certain distributions \cite{schwartz1978}. Thanks to this product, the solution of \eqref{eq:ode:intro} can be expressed through the $\star$-product inverse and its formulation as a sequence of integrals and differential equations; see \cite{Giscard2015,BonGis20,ProceedingsPaper2020,GiscardPozza2021,GiscardPozza2022}.
In \cite{PozVan_PANM21}, we illustrated that, by discretizing the $\star$-product with orthogonal functions, the solution of a scalar ODE is accessible by solving a linear system. 
In this work, we extend the results in \cite{PozVan_PANM21}, showing that, following the same principles, we can solve \eqref{eq:ode:intro} through a linear system. Moreover, we show that the linear system solution can be expressed as the solution of a matrix equation with a rank one right-hand side. Numerical experiments illustrate that the solution of the matrix equation can also be low-rank.

In Section \ref{sec:star}, we recall the $\star$-product definition and the related expression for the solution of an ODE. The expression is then discretized and approximated by the solution of a linear system. Section \ref{sec:mtxeq} shows how to transform the linear system into a matrix equation, and Section \ref{sec:conc} concludes the paper.


\section{Solution of an ODE by the $\star$-product}\label{sec:star}
We use the Heaviside theta function
\begin{equation*}
    \Theta(t-s) = \begin{cases}
                        1, \quad t \geq s \\
                        0, \quad t < s 
                \end{cases},
\end{equation*}
to rewrite \eqref{eq:ode:intro} in the following equivalent form
\begin{equation}\label{eq:ode:theta} 
\frac{d}{dt}U(t,s) = \tilde{A}(t) \Theta(t-s) U(t,s), \quad U(s,s)=I_N, \quad \text{ for } t,s\in \mathcal{I}.
\end{equation}
Note that $\Theta(t-s)$ endows the condition $t \geq s$ in equation  \eqref{eq:ode:theta} and that $U(t,s)$ is the bivariate function expressing the solutions of \eqref{eq:ode:intro} for every initial time $s \in \mathcal{I}$, with $U(t,s)=0$ for $t < s$.
From now on, we will denote with a tilde all the bivariate functions that are infinitely differentiable in both $t$ and $s$ over $\mathcal{I}$, i.e., $\tilde{f} \in C^\infty(\mathcal{I}\times \mathcal{I})$.
Moreover, we define the following class of functions
	\begin{equation*}
		\Sm(\mathcal{I}):= \left\{ f: f(t,s) = \tilde{f}(t,s) \Theta(t-s), \quad \tilde{f}\in C^{\infty}(\mathcal{I} \times \mathcal{I}) \right\}.
	\end{equation*}
Consider now the $N \times N$ matrices $A_1(t,s), A_2(t,s) \in (\Sm({\mathcal{I})})^{N \times N}$, i.e., matrices composed of elements from $\Sm(\mathcal{I})$. 
Then, the $\star$-product is defined as
\begin{equation}\label{eq:def:star}
  \big(A_2 \star A_1\big)(t,s) := \int_\mathcal{I} A_2(t,\tau) A_1(\tau, s) \, \text{d}\tau.
\end{equation}
The $\star$-product can be extended to a larger class of matrices composed of elements from the class $\D \supset \Sm(\mathcal{I})$, that is, the class of the superpositions of $\Theta(t-s)$, Dirac delta distribution $\delta(t-s)$, and Dirac delta derivatives described in \cite{ProceedingsPaper2020}. In such a class, $\delta(t-s) I_N$ is the $\star$-product identity, i.e., $A(t,s) \star  \delta(t-s)I_N = \delta(t-s)I_N \star A(t,s) = A(t,s)$.
  Moreover, in the larger class $\mathcal{D}(\mathcal{I})$, the $\star$-product admits inverses under certain conditions \cite{ProceedingsPaper2020}, i.e., for certain $f(t,s) \in \Sm$, there exists $f(t,s)^{-\star}$ such that $f(t,s)\star f(t,s)^{-\star} = f(t,s)^{-\star} \star f(t,s) = \delta(t-s)$.

Following \cite{Giscard2015}, the solution of \eqref{eq:ode:theta} can be expressed as
\begin{equation}\label{eq:star:sol}
    U(t,s) = \Theta(t-s) \star R_\star(A)(t,s),     
\end{equation}
where $A(t,s) = \tilde{A}(t) \Theta(t-s)$ and $R_\star(A)$ is the $\star$-resolvent of $A$, i.e.,
\begin{equation*}
    R_\star(A)(t,s) = \delta(t-s)I_N + \sum_{k=1}^\infty A(t,s)^{k \star},
\end{equation*}
with $A(t,s)^{k \star} = A\, \star \cdots \star \, A$, the $k$th power of the $\star$-product. Note that the series $\sum_{k=1}^\infty A(t,s)^{\star k}$ converges for every $A \in (\Sm(\mathcal{I}))^{N \times N}$.
Expression \eqref{eq:star:sol} hides an infinite series of nested integrals. However, as shown in \cite{PozVan_PANM21}, it is possible to approximate the $\star$-product by the usual matrix-matrix product in the scalar case. This approximation allows us to compute \eqref{eq:star:sol} more simply and cheaply. We recall its basics below. 

Without loss of generality, we set $\mathcal{I} = [0,1]$. Moreover, we consider the family of orthonormal shifted Legendre polynomials $\{ p_k \}_k$. Then, any $f(t,s) \in \Sm(\mathcal{I})$ can be expanded into the following series (e.g., \cite{LebSil72})
        \begin{equation}\label{eq:f:exp}
			f(t,s) = \sum_{k=0}^\infty \sum_{\ell=0}^\infty f_{k,\ell} \, p_k(t) p_\ell(s), \; t \neq s,\; t,s \in \mathcal{I}, \quad 
		    f_{k,\ell} = \int_\mathcal{I} \int_\mathcal{I} f(\tau,\rho) p_k(\tau) p_\ell(\rho) \; \textrm{d} \rho \; \textrm{d} \tau.
		\end{equation}
		By defining the \emph{coefficient matrix} $F_M$ and the vector $\phi_M(t)$ as
		\begin{equation}\label{eq:coeff:mtx}
		          F_M := \begin{bmatrix}
				f_{0,0} & f_{0,1} & \dots & f_{0,M-1}\\
				f_{1,0} & f_{1,1} & \dots & f_{1,M-1}\\
				\vdots & \vdots &  & \vdots\\
				f_{M-1,0} & f_{M-1,1} & \dots & f_{M-1,M-1}
		\end{bmatrix}, 
		\quad  \phi_M(t) :=
		 \begin{bmatrix}
			p_0(s)\\
			p_1(s)\\
			\vdots\\
			p_{M-1}(s)
		\end{bmatrix},
		\end{equation}
		the truncated expansion series can be written in the matrix form:
	\begin{align*}
	f_M(t,s) := \sum_{k=0}^{M-1} \sum_{\ell=0}^{M-1} f_{k,\ell} \, p_k(t) p_\ell(s) = \phi_M(t)^T F_M \, \phi_M(s).
		\end{align*}
  Let us consider the functions $f,g,h \in \Sm(\mathcal{I})$ so that $h = f \star g$, and the related coefficient matrices \eqref{eq:coeff:mtx}, respectively, $F_M,G_M,H_M$. Following \cite{ProceedingsPaper2020}, $H_M$ can be approximated by the expression
  \begin{equation}\label{eq:mtx:prod}
       H_M \approx \hat{H}_m := F_M G_M.
  \end{equation}
  
 Therefore, there is a connection between the $\star$-algebra over $\mathcal{D}(\mathcal{I})$ and the usual matrix algebra.
 The elements and operations which form the $\star$-algebra and the related elements and operations forming the usual matrix algebra are given in Table \ref{table:matrixOps} (in the first two columns for the scalar case); for more details, we refer to \cite{ProceedingsPaper2020}.

    The approximation in the scalar case can be easily extended to the matrix one. Indeed, if $A(t,s) = [a_{ij}(t,s)]_{i,j= 1}^N$ is an $N \times N$ matrix with elements $a_{ij}(t) \in \Sm(\mathcal{I})$, then for each $a_{ij}$, we can compute the related coefficient matrices $F_M^{(i,j)}$ \eqref{eq:coeff:mtx} obtaining the block matrix 
    \begin{equation}\label{eq:coeff:mtx:mtx}
         \mathcal{A}_M = \left[\begin{array}{ccccccc}
					\boxed{\begin{array}{c}
							F_M^{(1,1)}
					\end{array}} &  \dots & \boxed{\begin{array}{c}
							F_M^{(1,N)}
					\end{array}}\\
					\vdots         &   \ddots    & \vdots\\
					\boxed{\begin{array}{c}
							F_M^{(N,1)}
					\end{array}} & \dots & \boxed{\begin{array}{c}
							F_M^{(N,N)}
					\end{array}}\\
				\end{array}\right] \in \mathbb{C}^{MN \times MN}. 
\end{equation}
  Let us define the $N \times N$ matrices $A(t,s), B(t,s), C(t,s) \in (\Sm(\mathcal{I}))^{N \times N}$ so that $C(t,s) = A(t,s) \star B(t,s)$ and let their coefficient matrices \eqref{eq:coeff:mtx:mtx} be, respectively, $\mathcal{A}_M, \mathcal{B}_M, \mathcal{C}_M$. Then, analogously to the scalar case, $\mathcal{C}_M$ is approximated by 
  $$ \mathcal{C}_M \approx \hat{\mathcal{C}}_M := \mathcal{A}_M \mathcal{B}_M. $$
 As a consequence, also in the matrix case, the $\star$-algebra can be approximated by the usual matrix algebra, as summarized in the last two columns of Table \ref{table:matrixOps}.
 	\begin{table}[!ht]
		\begin{tabular}{ll|ll}
			$f(t,s) \in \Sm(\mathcal{I})$ & $F_M \in \mathbb{C}^{M \times M}$ & $A(t,s) \in (\Sm(\mathcal{I}))^{N \times N}$ & $\mathcal{A}_M \in \mathbb{C}^{MN \times MN}$ \\ \hline
			$\star$-operation/elements & matrix operation/elements & 	$\star$-operation/elements & matrix operation/elements\\				\hline
			$q = f \star g$             &  $Q_M = F_M G_M$ & $C = A \star B$             &  $\mathcal{C}_M = \mathcal{A}_M \mathcal{B}_M$ \\
			$f + g$           & $F_M+G_M$  & $A + B$           & $\mathcal{A}_M+\mathcal{B}_M$  \\
			$1_{\star} := \delta(t-s)$ &   $I_M$, identity matrix       & $1_{\star} := \delta(t-s) I_N$ &   $I_{MN}$, identity matrix       \\
			$f^{\star-1}$                    & $F_M^{-1}$  &  $A^{\star-1}$                    & $\mathcal{A}_M^{-1}$    \\
			$R_\star(f):=(1_{\star}- f)^{\star-1}$   &  $R(F_M) := (I_M-F_M)^{-1}$ & $R_\star(f):=(1_{\star}- A)^{\star-1}$   &  $R(\mathcal{A}_M) := (I_{MN}-\mathcal{A}_M)^{-1}$ 
		\end{tabular}
		\caption{The $\star$-algebra operations and the corresponding matrix algebra operation after discretization, scalar case (first two columns), matrix case (last two columns).}
		\label{table:matrixOps}
	\end{table}
 
   The matrix-valued function $U(t,s)$ in \eqref{eq:ode:theta} is composed of elements from $\Sm(\mathcal{I})$. Therefore, we can define the related coefficient matrix $\mathcal{U}_M$ as in \eqref{eq:coeff:mtx:mtx}. Then, expression~\eqref{eq:star:sol} can be approximated by
   \begin{equation*}
       \mathcal{U}_M \approx (I_N \otimes T_M) (I_{MN}-\mathcal{A}_M)^{-1},
   \end{equation*}
   where $\otimes$ is the Kronecker product, $T_M$ is the coefficient matrix of $\Theta(t-s)$, and $\mathcal{A}_M$ is the coefficient matrix of $\tilde{A}(t)\Theta(t-s)$, with $\tilde{A}(t)$ from \eqref{eq:ode:theta}.
   Moreover, we can approximate the solution of \eqref{eq:ode:theta} for $s=0$ by the formula:
   \begin{align*}
       U(t,0) \approx \phi_M(t)^T \mathcal{U}_M \, \phi_M(0) &= (I_N \otimes \phi_M(t)^T) (I_N \otimes T_M) (I_{MN}-\mathcal{A}_M)^{-1}  (I_N \otimes \phi_M(0)) \\
                                                   &= (I_N \otimes \phi_M(t)^T T_M) (I_{MN}-\mathcal{A}_M)^{-1}  (I_N \otimes \phi_M(0)).
   \end{align*}
   Note that, as explained in \cite{PozVan_PANM21}, the approximation converges quickly enough to the solution only when $s$ is the left endpoint of the interval $\mathcal{I}$, i.e., $s=0$.
   
   In practical situations, the initial time $s$ of the evolution is fixed ($s=0$), and the initial condition is given as a vector $v \in \mathbb{C}^N$.
   Then, we get the simpler problem,
   \begin{equation}\label{eq:ode:vect}
\frac{d}{dt}u(t) = \tilde{A}(t) \Theta(t-s) u(t), \quad u(0)=v, \quad \text{ for } t,s\in \mathcal{I},
\end{equation}
where the solution $u(t)$ is an $N$-size vector. Thus, $u(t)$ is approximated by:
   \begin{align*}
       u(t) &\approx (I_N \otimes \phi_M(t)^T T_M) (I_{MN}-\mathcal{A}_M)^{-1} (I_N \otimes \phi_M(0))  \, v  \\
            &\approx (I_N \otimes \phi_M(t)^T T_M) (I_{MN}-\mathcal{A}_M)^{-1} (v \otimes \phi_M(0)).
   \end{align*}
   Then, solving the linear system
     \begin{equation}\label{eq:linesys}
         (I_{MN}-\mathcal{A}_M) x = v \otimes \phi_M(0),
     \end{equation}
   one can approximate the solution of \eqref{eq:ode:vect} in terms of its expansion coefficients $u_M := (I_N \otimes T_M) x$, that is,
   \begin{equation}\label{eq:approx}
        u(t) \approx \hat{u}(t):= (I_N \otimes \phi_M(t)^T) u_M.
   \end{equation}

\subsection{Numerical examples}\label{sec:num:ex}
Given a random vector $v$ with elements in $[0,1]$, we aim to compute the bilinear form $v^T u(t)$ obtained by solving the following ODE system
\begin{equation}\label{eq:ode:MAS}
    \frac{d}{dt}u(t) = - 2\sqrt{-1}\pi \tilde{H}(t) u(t), \quad u(0)=v, \quad \text{ for } t \in [0, T].
\end{equation}
This system of ODEs comes from Experiment 2 (Strong coupling) in \cite{CiPetal22}, and $v^T u(t)$ 
 represents an NMR experiment with a magic angle spinning (MAS) for $k$ spins; see, e.g., \cite{HaSp98}.
The so-called Hamiltonian $\tilde{H}(t)$ is a $2^k \times 2^k$ matrix-valued function and has the form
\begin{equation}\label{eq:hamiltonian}
    \tilde{H}(t) = D + B(\cos(2\pi\nu t) + \cos(4\pi\nu t)),
\end{equation}
with $D, B$ sparse matrices described in \cite{CiPetal22}.
In our experiments, we set $T = 10^{-3}$, $\nu=10^4$, and $k=4, 7, 10$, so obtaining three systems with exponentially increasing sizes.

\begin{figure}[ht]
\includegraphics[width=5.4cm]{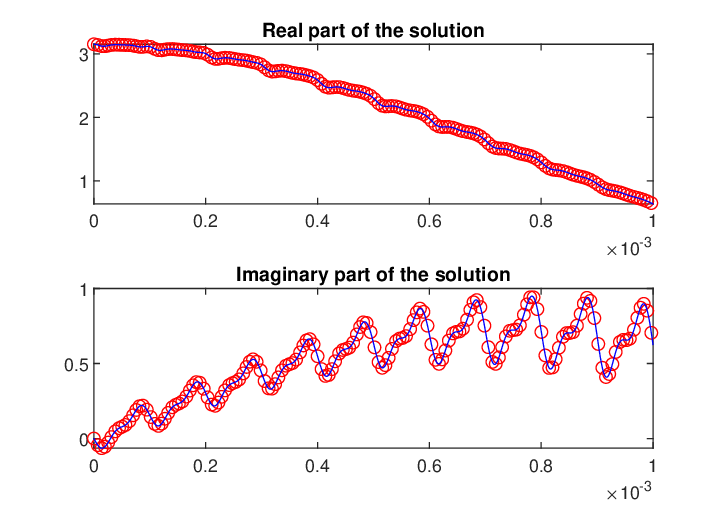}
\includegraphics[width=5.4cm]{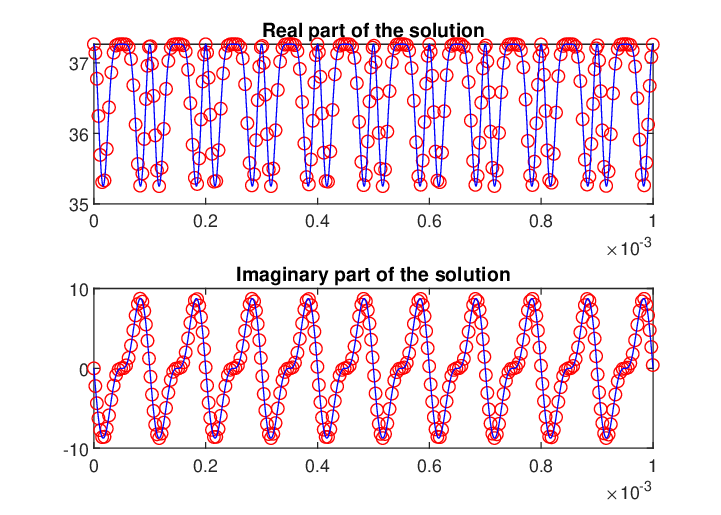}
\includegraphics[width=5.4cm]{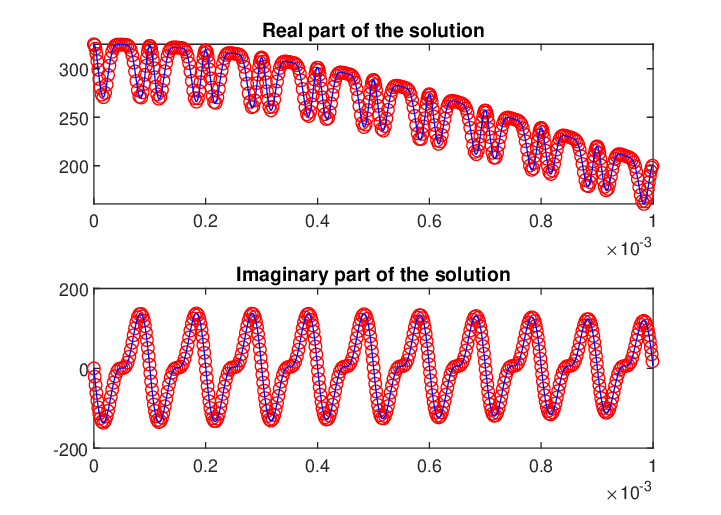}
\caption{Real and imaginary parts of $v^T u(t)$ approximations, with $u(t)$ the solution of \eqref{eq:ode:MAS}. The red circles represent approximation $v^T \hat{u}(t)$ from \eqref{eq:approx}, while the blue line represents the \texttt{ode45} approximation. From left to right, $k=4, 7, 10$.}
\label{fig:sol}
\end{figure}

\begin{figure}[ht]
\includegraphics[width=5.4cm]{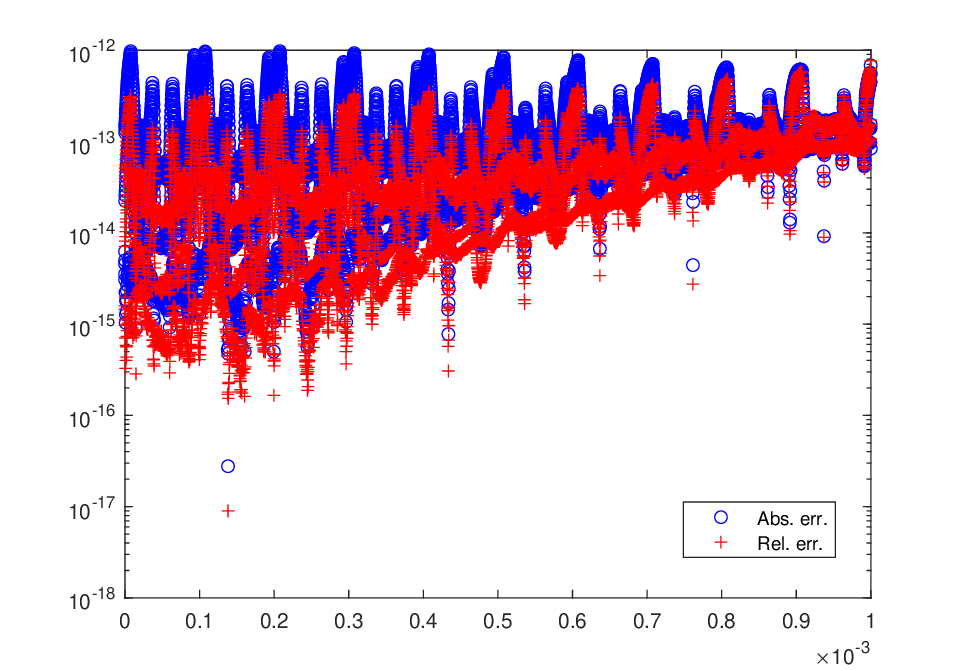}
\includegraphics[width=5.4cm]{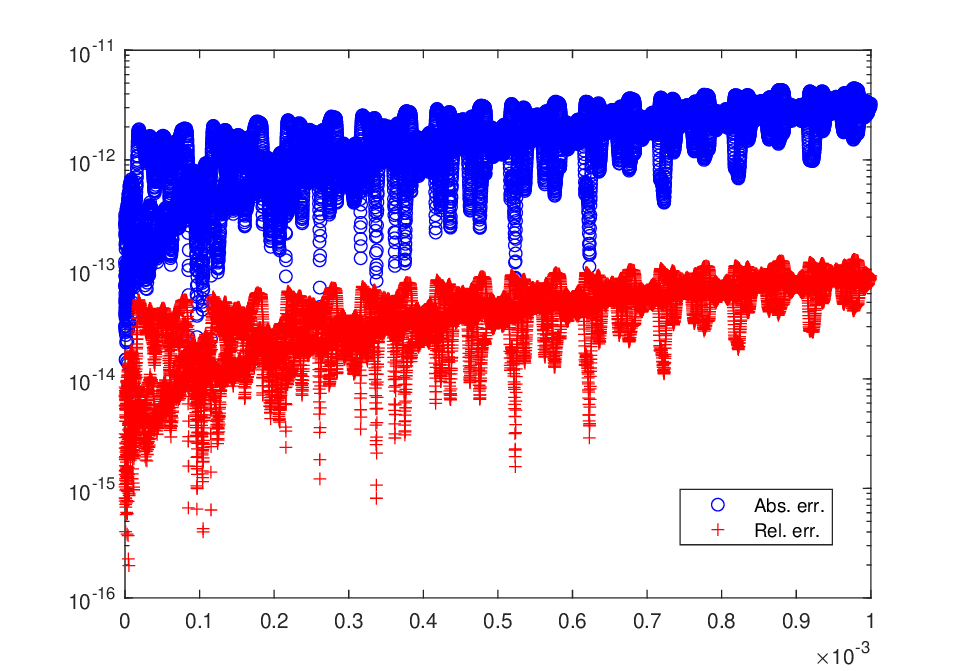}
\includegraphics[width=5.4cm]{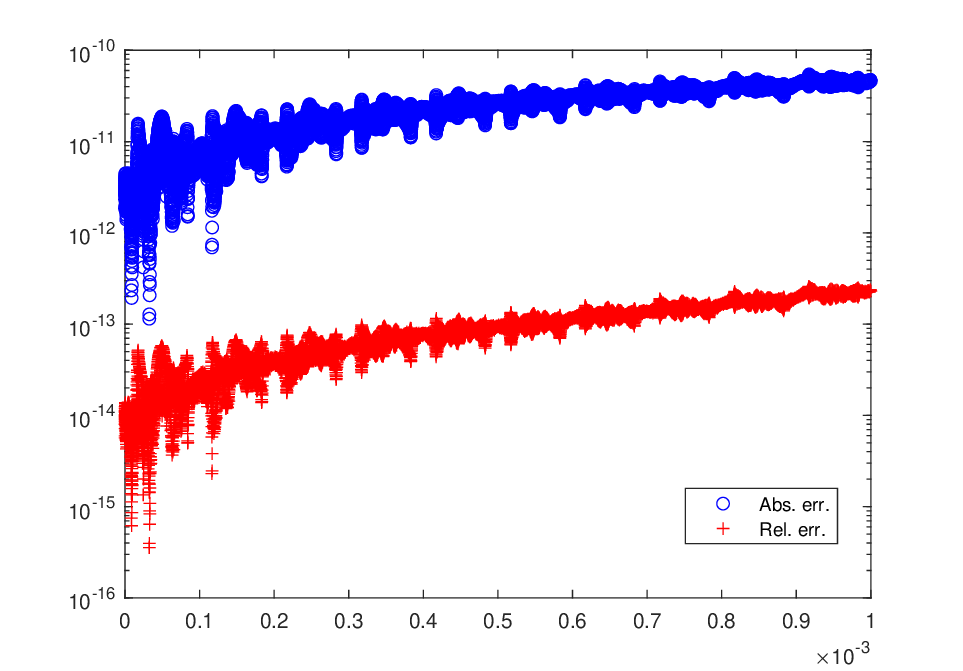}
\caption{Absolute (blue circles) and relative (red crosses) errors of approximation $v^T \hat{u}(t)$, with $\hat{u}(t)$ from \eqref{eq:approx}. From left to right, $k=4, 7, 10$.}
\label{fig:err}
\end{figure}

The approximated solution $\hat{u}(t)$ \eqref{eq:approx} is computed by solving the linear system \eqref{eq:linesys}\footnote{The matrices $F^{(i,j)}_M$ in the block coefficient matrix \eqref{eq:coeff:mtx:mtx} are numerically banded with bandwidth $b_{i,j}$. In order to avoid error accumulation, the last $b_{i,j}$ rows of each $F^{(i,j)}_M$ have been set equal to zero; see \cite{PozVan_PANM21}.} with $M= 1000$.
The numerical experiments were performed using MatLab R2022a, and the linear systems were solved by the MatLab GMRES method implementation, \texttt{gmres}, with tolerance set to $1e-15$.
In Figure~\ref{fig:sol}, we compare the approximated bilinear form $v^T \hat{u}(t)$ with the solution obtained by the MatLab function \texttt{ode45} with relative and absolute tolerance set to $3e-14$.
Figure~\ref{fig:err} reports the corresponding relative and absolute errors over the interval $[0,T]$ (the reference for the error is again the \texttt{ode45} solution).
In all the experiments, GMRES stopped after a maximum of 27 iterations (for the cases $k=7,10$ due to residual stagnation).
The numerical results show that the method is able to compute the solution with accuracy comparable with a well-established method.


\section{Matrix equation formulation}\label{sec:mtxeq}
\begin{figure}[ht]
\includegraphics[width=5.4cm]{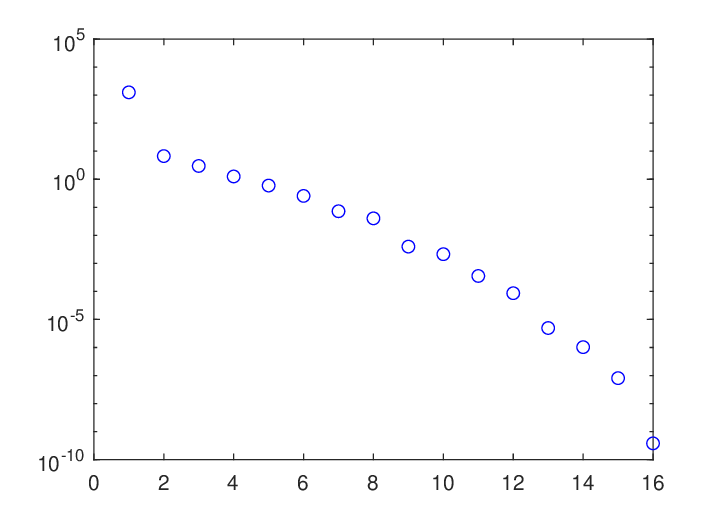}
\includegraphics[width=5.4cm]{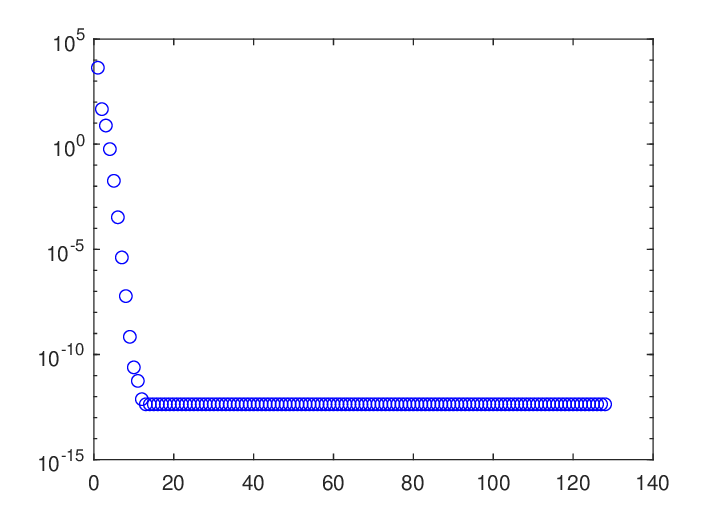}
\includegraphics[width=5.4cm]{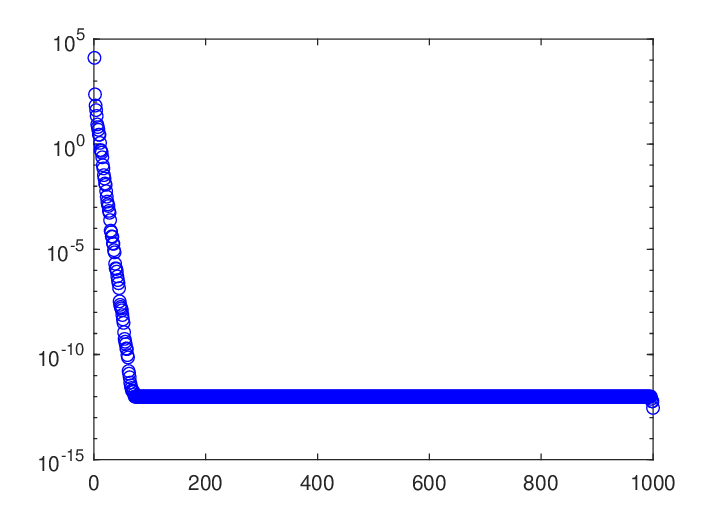}
\caption{Singular values of the matrix $X$, with $x = vec(X)$ the solution of \eqref{eq:linesys} for the examples in Section \ref{sec:num:ex}. From left to right, $k=4,7,10$.}
\label{fig:svd}
\end{figure}

The matrix-valued function $\tilde{A}(t)$ in \eqref{eq:ode:intro} can always be written in the form
\begin{equation}\label{eq:Asum}
      \tilde{A}(t) = \sum_{k=1}^d A_k \tilde{f}_k(t),
\end{equation}
 with $\tilde{f}_1, \dots, \tilde{f}_d$ distinct scalar functions and $A_1, \dots, A_d$ constant matrices.
In many applications, $d$ is small. For instance, in the examples from Section \ref{sec:num:ex}, we have $d=2$.
Then, exploiting expression \eqref{eq:Asum}, the (block) coefficient matrix \eqref{eq:coeff:mtx:mtx} of $A(t,s)= \tilde{A}(t)\Theta(t-s)$ becomes
\begin{equation*}
      \mathcal{A}_M = \sum_{k=1}^d A_k \otimes F_M^{(k)},
\end{equation*}
with $F_M^{(k)}$ the coefficient matrix \eqref{eq:coeff:mtx} of $\tilde{f}_k(t)$. The solution $x$ of the linear system \eqref{eq:linesys} can, hence, 
be rewritten in terms of the solution $X$ of the following matrix equation
\begin{equation}\label{eq:mtxeq}
       X - \sum_{k=1}^d F_M^{(k)} X A_k^T = \phi_M(0) b^T, \quad x = vec(X),
\end{equation}
where  $vec(X)$ denotes the vectorization of $X$, i.e., the vector obtained by stacking the columns of $X$ into a single vector. 
The matrix equation \eqref{eq:mtxeq} has a rank $1$ right-hand side $\phi_M(0) b^T$. 
This suggests that the solution $X$ may have a low numerical rank. 
Figure \ref{fig:svd} reports the computed singular values of $X$, where $x = vec(X)$ is the linear system solution of each of the experiments performed in Section \ref{sec:num:ex}.
For $k=4$, the solution $X$ is full rank, while for $k=7, 10$, the numerical rank of $X$ is, respectively, $12, 72$ 
(we consider as numerical rank the index of the last singular value before the stagnation visible in the plots).
Clearly, this preliminary study shows that the numerical rank of $X$ increases slowly with the size of $X$.

\section{Discussions and conclusion}\label{sec:conc}
In this work, we present a new method for solving systems of  non-autonomous linear ODEs. The method is based on the solution of a linear system that can be rewritten as a matrix equation.
Several examples illustrate that the method is able to compute the solution with accuracy comparable to the well-established Runge-Kutta method implemented by the MatLab function \texttt{ode45}. 
Moreover, the experiments show that the solution of the matrix equation is a numerical low-rank matrix when the ODE system is large enough. This may be exploited using projection methods with low-rank techniques (see, e.g., \cite{Sim07,KurPal21}). In \cite{CiPetal22}, we also show that matrix $\mathcal{A}_M$ in \eqref{eq:coeff:mtx:mtx} can be compressed by the Tensor Train decomposition (note that \cite{CiPetal22} uses a different family of orthogonal functions instead of the Legendre polynomials). A Tensor Train approach may further reduce the memory and computational cost of the method.
Another possible approach could be extrapolation methods able to exploit the dependence of equation~\eqref{eq:linesys} on $s$; see, e.g., \cite{BreRed91,BuoKarPoz15}.

Overall, the results suggest that the presented method may be an effective solver for large-to-huge systems of ODEs once we are able to exploit the solution's low-rank structure and the other mentioned properties. We are currently investigating these possible approaches.

%

\vspace{\baselineskip}

\bibliographystyle{pamm}
\bibliography{biblio}

\end{document}